\magnification=1200

\hsize=11.25cm    
\vsize=18cm     
\parindent=12pt   \parskip=5pt     

\hoffset=.5cm   
\voffset=.8cm   

\pretolerance=500 \tolerance=1000  \brokenpenalty=5000

\catcode`\@=11

\font\eightrm=cmr8         \font\eighti=cmmi8
\font\eightsy=cmsy8        \font\eightbf=cmbx8
\font\eighttt=cmtt8        \font\eightit=cmti8
\font\eightsl=cmsl8        \font\sixrm=cmr6
\font\sixi=cmmi6           \font\sixsy=cmsy6
\font\sixbf=cmbx6

\font\tengoth=eufm10 
\font\eightgoth=eufm8  
\font\sevengoth=eufm7      
\font\sixgoth=eufm6        \font\fivegoth=eufm5

\skewchar\eighti='177 \skewchar\sixi='177
\skewchar\eightsy='60 \skewchar\sixsy='60

\newfam\gothfam           \newfam\bboardfam

\def\tenpoint{
  \textfont0=\tenrm \scriptfont0=\sevenrm \scriptscriptfont0=\fiverm
  \def\rm{\fam\z@\tenrm}
  \textfont1=\teni  \scriptfont1=\seveni  \scriptscriptfont1=\fivei
  \def\oldstyle{\fam\@ne\teni}\let\old=\oldstyle
  \textfont2=\tensy \scriptfont2=\sevensy \scriptscriptfont2=\fivesy
  \textfont\gothfam=\tengoth \scriptfont\gothfam=\sevengoth
  \scriptscriptfont\gothfam=\fivegoth
  \def\goth{\fam\gothfam\tengoth}
  
  \textfont\itfam=\tenit
  \def\it{\fam\itfam\tenit}
  \textfont\slfam=\tensl
  \def\sl{\fam\slfam\tensl}
  \textfont\bffam=\tenbf \scriptfont\bffam=\sevenbf
  \scriptscriptfont\bffam=\fivebf
  \def\bf{\fam\bffam\tenbf}
  \textfont\ttfam=\tentt
  \def\tt{\fam\ttfam\tentt}
  \abovedisplayskip=12pt plus 3pt minus 9pt
  \belowdisplayskip=\abovedisplayskip
  \abovedisplayshortskip=0pt plus 3pt
  \belowdisplayshortskip=4pt plus 3pt 
  \smallskipamount=3pt plus 1pt minus 1pt
  \medskipamount=6pt plus 2pt minus 2pt
  \bigskipamount=12pt plus 4pt minus 4pt
  \normalbaselineskip=12pt
  \setbox\strutbox=\hbox{\vrule height8.5pt depth3.5pt width0pt}
  \let\bigf@nt=\tenrm       \let\smallf@nt=\sevenrm
  \normalbaselines\rm}

\def\eightpoint{
  \textfont0=\eightrm \scriptfont0=\sixrm \scriptscriptfont0=\fiverm
  \def\rm{\fam\z@\eightrm}
  \textfont1=\eighti  \scriptfont1=\sixi  \scriptscriptfont1=\fivei
  \def\oldstyle{\fam\@ne\eighti}\let\old=\oldstyle
  \textfont2=\eightsy \scriptfont2=\sixsy \scriptscriptfont2=\fivesy
  \textfont\gothfam=\eightgoth \scriptfont\gothfam=\sixgoth
  \scriptscriptfont\gothfam=\fivegoth
  \def\goth{\fam\gothfam\eightgoth}
  
  \textfont\itfam=\eightit
  \def\it{\fam\itfam\eightit}
  \textfont\slfam=\eightsl
  \def\sl{\fam\slfam\eightsl}
  \textfont\bffam=\eightbf \scriptfont\bffam=\sixbf
  \scriptscriptfont\bffam=\fivebf
  \def\bf{\fam\bffam\eightbf}
  \textfont\ttfam=\eighttt
  \def\tt{\fam\ttfam\eighttt}
  \abovedisplayskip=9pt plus 3pt minus 9pt
  \belowdisplayskip=\abovedisplayskip
  \abovedisplayshortskip=0pt plus 3pt
  \belowdisplayshortskip=3pt plus 3pt 
  \smallskipamount=2pt plus 1pt minus 1pt
  \medskipamount=4pt plus 2pt minus 1pt
  \bigskipamount=9pt plus 3pt minus 3pt
  \normalbaselineskip=9pt
  \setbox\strutbox=\hbox{\vrule height7pt depth2pt width0pt}
  \let\bigf@nt=\eightrm     \let\smallf@nt=\sixrm
  \normalbaselines\rm}

\tenpoint

\def\pc#1{\bigf@nt#1\smallf@nt}         \def\pd#1 {{\pc#1} }

\frenchspacing

\def\raggedbottom{\topskip 10pt plus 36pt\r@ggedbottomtrue}

\def\pointir{\unskip . --- \ignorespaces}

\def\Medbreak{\vskip-\lastskip\medbreak}

\long\def\th#1 #2\enonce#3\endth{
   \Medbreak\noindent
   {\pc#1} {#2\unskip}\pointir{\it #3}\smallskip}

\def\decale#1{\smallbreak\hskip 28pt\llap{#1}\kern 5pt}
\def\decaledecale#1{\smallbreak\hskip 34pt\llap{#1}\kern 5pt}
\def\puce{\smallbreak\hskip 6pt{$\scriptstyle\bullet$}\kern 5pt}

\def\eqalign#1{\null\,\vcenter{\openup\jot\m@th\ialign{
\strut\hfil$\displaystyle{##}$&$\displaystyle{{}##}$\hfil
&&\quad\strut\hfil$\displaystyle{##}$&$\displaystyle{{}##}$\hfil
\crcr#1\crcr}}\,}

\catcode`\@=12

\showboxbreadth=-1  \showboxdepth=-1

\newcount\numerodesection \numerodesection=1
\def\section#1{\bigbreak
 {\bf\number\numerodesection.\ \ #1}\nobreak\medskip
 \advance\numerodesection by1}

\mathcode`A="7041 \mathcode`B="7042 \mathcode`C="7043 \mathcode`D="7044
\mathcode`E="7045 \mathcode`F="7046 \mathcode`G="7047 \mathcode`H="7048
\mathcode`I="7049 \mathcode`J="704A \mathcode`K="704B \mathcode`L="704C
\mathcode`M="704D \mathcode`N="704E \mathcode`O="704F \mathcode`P="7050
\mathcode`Q="7051 \mathcode`R="7052 \mathcode`S="7053 \mathcode`T="7054
\mathcode`U="7055 \mathcode`V="7056 \mathcode`W="7057 \mathcode`X="7058
\mathcode`Y="7059 \mathcode`Z="705A


\def\diagram#1{\def\normalbaselines{\baselineskip=0pt\lineskip=5pt}
\matrix{#1}}

\def\vfl#1#2#3{\llap{$\textstyle #1$}
\left\downarrow\vbox to#3{}\right.\rlap{$\textstyle #2$}}

\def\hfl#1#2#3{\smash{\mathop{\hbox to#3{\rightarrowfill}}\limits
^{\textstyle#1}_{\textstyle#2}}}

\def\ogoth{{\goth o}}

\def\pgoth{{\goth p}}

\def\Q{{\bf Q}}

\def\R{{\bf R}}

\def\N{{\bf N}}

\def\Z{{\bf Z}}

\def\F{{\bf F}}
\def\Fp{\F_p}

\def\Card{\mathop{\rm Card}\nolimits}
\def\Gal{\mathop{\rm Gal}\nolimits}

\def\to{\rightarrow}

\def\droite#1{\,\hfl{#1}{}{8mm}\,}

\def\mod{\mathop{\rm mod.}\nolimits}
\def\pmod#1{\;(\mod#1)}

\def\boxit#1{\vbox{\hrule\hbox{\vrule\kern1pt
       \vbox{\kern1pt#1\kern1pt}\kern1pt\vrule}\hrule}}
\def\cqfd{\hfill\boxit{\phantom{\i}}}

\newcount\numero \numero=1
\def\numeroter{{({\oldstyle\number\numero})}\ \advance\numero by1}

\newcount\refno 
\long\def\ref#1:#2<#3>{                                        
\global\advance\refno by1\par\noindent                              
\llap{[{\bf\number\refno}]\ }{#1} \pointir{\it #2} #3\goodbreak }

\newcount\refno 
\long\def\ref#1:#2<#3>{                                        
\global\advance\refno by1\par\noindent                              
\llap{[{\bf\number\refno}]\ }{#1} \pointir{\it #2} #3\goodbreak }

\def\citer#1(#2){[{\bf\number#1}\if#2\empty\relax\else,\ {#2}\fi]}

\def\boxit#1{\vbox{\hrule\hbox{\vrule\kern1pt
       \vbox{\kern1pt#1\kern1pt}\kern1pt\vrule}\hrule}}
\def\cqfd{\hfill\boxit{\phantom{\i}}}

\newbox\bibbox
\setbox\bibbox\vbox{\bigbreak
\centerline{{\pc BIBLIOGRAPHY}}
\nobreak

\ref{\pc BOYLAN} (Hatice):
Question 248665,
<mathoverflow.net/q/248665.>
\newcount\boylan \global\boylan=\refno

\ref{\pc DALAWAT} (Chandan):
Local discriminants, kummerian extensions, and elliptic curves,
<J.\ Ramanujan math.\ soc.\ {\bf 25} (2010) 1, 25--80.  arXiv\string:0711.3878.>      
\newcount\locdisc \global\locdisc=\refno

\ref{\pc DALAWAT} (Chandan):
Further remarks on local discriminants, 
<J.\ Ramanujan math.\ soc.\ {\bf 25} (2010) 4, 397--417.  arXiv\string:0909.2541.>    
\newcount\further \global\further=\refno

\ref{\pc SERRE} (Jen-Pierre):
Corps locaux,
<Publications de l'Universit{\'e} de Nancago {\sevenrm VIII}, Hermann,
Paris, 1968, 245 p.>
\newcount\corpslocaux \global\corpslocaux=\refno

} 

\centerline{\bf Reciprocity and orthogonality} 
\bigskip\bigskip 
\centerline{Chandan Singh Dalawat} 
\centerline{Harish-Chandra Research Institute}
\centerline{Chhatnag Road, Jhunsi, Allahabad 211019, India} 
\centerline{\tt dalawat@gmail.com}

\bigskip\bigskip

{{\bf Abstract}.  Let $p$ be a prime and let $K$ be a finite extension
  of the field $\Q_p$ of $p$-adic numbers such that the group
  ${}_pK^\times$ has order~$p$. The $\F_p$-space $K^\times\!/K^{\times
  p}$ carries a natural filtration coming from the valuation on $K$,
  and a natural bilinear pairing coming from the reciprocity
  isomorphism for the exponent~$p$.  We determine the orthogonal
  filtration for this pairing.  We also prove the analogous result for
  $p$-fields of characteristic~$p$.
 
\footnote{}{{\it MSC2010~:} Primary  11S15, 11S31} 
\footnote{}{{\it Keywords~:} Hilbertian pairing, Reciprocity
  isomorphism, Orthogonality}
}

 \bigskip\bigskip


 \numeroter Let $p$ be a prime number and let $K$ be a $p$-field --- a
 field complete for a discrete valuation with finite residue field of
 characteristic~$p$.  Let $M$ be the maximal abelian extension of $K$
 of exponent $p$, $G=\Gal(M|K)$,
 $\overline{K^\times}=K^\times\!/K^{\times p}$, and
$\rho:\overline{K^\times}\to G$ the reciprocity isomorphism (for the
 exponent~$p$). 
 If $K$ has characteristic~$0$ and ${}_pK^\times$ has
 order~$p$, then there is a natural pairing
 $$
 G\times\overline{K^\times}\to{}_pK^\times,\quad
 (\sigma,\bar x)\mapsto\sigma(y)y^{-1}\ \  (y\in M, y^p=x).
 $$
  Similarly, if
 $K$ has characteristic~$p$, then there is a natural pairing
 $$
 G\times\overline{K^+}\to\F_p,\quad
 (\sigma,\bar x)\mapsto\sigma(y)-y\ \  (y\in M, \wp(y)=x),
 $$
where $\overline{K^+}=K^+\!/\wp(K^+)$ and $\wp(z)=z^p-z$ for $z$ in
 any $\F_p$-algebra.  

\numeroter When combined with the reciprocity
isomorphism~$\rho:\overline{K^\times}\to G$, we get the {\it
  hilbertian pairing\/}
$$
\overline{K^\times}\times\overline{K^\times}\to{}_pK^\times\ (\hbox{resp.\ }
\overline{K^\times}\times\overline{K^+}\to\F_p).
$$
The main aim of this Note is to show that the filtrations on the two
factors, coming from the filtration(s) on $K^\times$ (resp.~on
$K^\times$ and $K^+$), are orthogonal to each other for this pairing.
This will be soon made precise.

\bigbreak
{\bf 1.  Notation}
\bigskip

\numeroter ($p$, $K$, $\ogoth$, $\pgoth$, $k$).  Throughout, $p$ is a
prime number and $K$ is a $p$-field --- a field complete for a
discrete valuation with finite residue field of characteristic~$p$.
The ring of integers of $K$ is $\ogoth$, the unique maximal ideal of
$\ogoth$ is $\pgoth$, and the residue field is $k=\ogoth/\pgoth$.

\numeroter ($e$, $f$).  The (absolute) residual degree of $K$ is
$f=[k:\F_p]$.  If $K$ is an extension of $\Q_p$ , then
$e=[K:\Q_p]f^{-1}$ is the (absolute) ramification index~; if $K$ has
characteristic~$p$, then $e=+\infty$.  We often say ``\thinspace
$e<+\infty$\thinspace'' to mean that $K$ has characteristic~$0$.

\numeroter ($\overline{K^\times}$, $\bar U_i$).  For every $i>0$, let
$U_i=1+\pgoth^i$ be the kernel of
$\ogoth^\times\to(\ogoth/\pgoth^i)^\times$ and denote the image of
$U_i$ in the $\F_p$-space $\overline{K^\times}=K^\times\!/K^{\times
  p}$ by $\bar U_i$.  Note that the image 
$\ogoth^\times\!/\ogoth^{\times p}$ of $\ogoth^\times$ in
$\overline{K^\times}$ is equal to $\bar U_1$ because
$\ogoth^\times=U_1.k^\times$ and $k^{\times p}=k^\times$.  The image
of $x\in K^\times$ in $\overline{K^\times}$ is denoted $\bar x$.

\numeroter ($\wp$, $\overline{K^+}$, $\overline{\pgoth^i}$,
$\bar\ogoth$).  For any $\F_p$-algebra $A$, denote by $\wp:A\to A$ the
endomorphism $\wp(x)=x^p-x$ of the additive group $A^+$. If
$e=+\infty$ (so that $K$ is an $\F_p$-algebra), we put
$\overline{K^+}=K^+/\wp(K^+)$ and, for every $i\in\Z$, denote the
image of $\pgoth^i$ in the $\F_p$-space $\overline{K^+}$ by
$\overline{\pgoth^i}$.  Sometimes we also denote $\overline{\pgoth^0}$
by $\bar\ogoth$~; it equals $\ogoth^+\!/\wp(\ogoth^+)$.  The image
of $x\in K^+$ in $\overline{K^+}$ is denoted $\bar x$.

\numeroter ($b_p^{(i)}$).  We denote by $b_p^{(i)}$ ($i>0$) be the
sequence of positive integers $\not\equiv0\pmod p$, namely
$b_p^{(i)}=i+\lfloor(i-1)/(p-1)\rfloor$, where, for every $x\in\R$,
$\lfloor x\rfloor$ is the largest integer in the interval $]-\infty,x]$.

\numeroter ($c$).  If $e<+\infty$ and $e\equiv0\pmod{p-1}$ (for
example when ${}_pK^\times$ has order~$p$), then we abbreviate
$c=(p-1)^{-1}e$.  Note that in this case $pc=b_p^{(e)}+1=e+c$.

\numeroter ($v$, $\bar v$).  The surjective valuation $K^\times\to\Z$
is denoted by $v$ (so that $e=v(p)$)~; it induces as isomorphism $\bar
v:\overline{K^\times}/\bar U_1\to\Z/p\Z$.  For this reason, in order
to determine the structure the the filtered $\F_p$-space
$\overline{K^\times}$, it is enough to study the filtered $\F_p$-space
$\bar U_1$.

\bigbreak
{\bf 2.  The filtered $\F_p$-space $\bar U_1$}
\bigskip
    
\numeroter Let us determine the filtration on the $\F_p$-space $\bar
U_1$.  We will see that the dimension $d=\dim_{\F_p}\bar U_1$ is
finite or infinite according as $e<+\infty$ or $e=+\infty$. When
$e<+\infty$, one has $d=[K:\Q_p]+\dim_{\F_p}({}_pK^\times)$, so there
are two subcases according as ${}_pK^\times$ is trivial or has order~$p$.
More precisely, by studying what the endomorphism $(\ )^p:U_1\to U_1$
does to the filtration on $U_1$, one determines the filtration on
$\bar U_1$ as follows.

\numeroter {\it Suppose that\/ $e=+\infty$.  For\/ $i>0$, the
  inclusion\/ $\bar U_{i+1}\subset\bar U_i$ is an equality if\/
  $i\equiv0\pmod p$, and has codimension~$f$ if\/ $i\not\equiv0\pmod p$.}
\cqfd 

\numeroter {\it Suppose that\/ $e<+\infty$ and\/ ${}_pK^\times$ is
  trivial.  Then\/ $\bar U_i$ is trivial for\/ $i>b_p^{(e)}$.  For\/
  $i\in[1,b_p^{(e)}]$, the inclusion\/ $\bar U_{i+1}\subset\bar U_i$
  is an equality if\/ $i\equiv0\pmod p$, and has codimension~$f$ if\/
  $i\not\equiv0\pmod p$.}
\cqfd 

\numeroter {\it Suppose that\/ $e<+\infty$ and\/ ${}_pK^\times$ has
  order\/~$p$.  Then\/ $\bar U_i$ is trivial for\/ $i>pc$ and\/ $\bar
  U_{pc}$ has order~$p$.  For\/ $i\in[1,pc[$, the inclusion\/ $\bar
  U_{i+1}\subset\bar U_i$ is an equality if\/ $i\equiv0\pmod p$, and
  has codimension~$f$ if\/ $i\not\equiv0\pmod p$.}  \cqfd

\numeroter {\it Remark}.  Whenever $\bar U_{i+1}$ has
  codimension~$f$ in $\bar U_i$, the quotient $\bar U_i/\bar U_{i+1}$
  is canonically isomorphic to $U_i/U_{i+1}$.

\numeroter Denote the reduction map $\ogoth\to k$ by $a\mapsto\hat a$
and let $S:k\to\F_p$ be the trace map.  If $e<+\infty$ and if
$\zeta\in{}_pK^\times$ has order~$p$, then $v(p\pi)=pc$, where
$\pi=1-\zeta$, and the map $\overline{1+ap\pi}\mapsto\zeta^{S(\hat
a)}$ ($a\in\ogoth$) is an isomorphism $\bar U_{pc}\to{}_pK^\times$,
independent of $\zeta$.  For proofs and more information,
see \citer\locdisc(Prop.~42), for example.

\bigbreak
{\bf 3.  The filtered $\F_p$-space $\overline{K^+}$}
\bigskip
    
\numeroter  Suppose that $e=+\infty$.  In analogy with the foregoing,
by studying what the endomorphism $\wp:K^+\to K^+$ does to the
filtration on $K^+$ (by the powers $\pgoth^i$ ($i\in\Z$) of $\pgoth$),
one determines the filtration on $\overline{K^+}$.

\numeroter {\it For every\/ $i>0$, one has\/ $\overline{\pgoth^i}=\{0\}$.
The group $\overline{\pgoth^0}=\bar\ogoth$ has order~$p$.  For every
$i<0$, the inclusion\/
$\overline{\pgoth^{i+1}}\subset\overline{\pgoth^i}$ is an equality
if\/ $i\equiv0\pmod p$, and has codimension~$f$ if\/
$i\not\equiv0\pmod p$.}
\cqfd

\numeroter  {\it Remark}.  Denote the passage to the quotient
$\ogoth\to\bar\ogoth$ (resp.~the reduction map $\ogoth\to k$) by
$x\mapsto\bar x$ (resp.~$x\mapsto\hat x$), and the trace map
$k\to\F_p$ by $S$.  Then the map $\bar a\mapsto S(\hat a)$
($a\in\ogoth$) is an isomorphism $\bar\ogoth\to\F_p$.  For proofs and
more information, see \citer\further(Prop.~11), for example.
  
\bigbreak
{\bf 4.  Breaks and levels}
\bigskip

\numeroter Let $E$ be a cyclic extension of $K$ of
degree~$p$.  The ramification filtration on the group $G=\Gal(E|K)$
has a unique {\it break\/} $\varepsilon(E)$ (the integer $j$ such that
$G^j=G$, $G^{j+1}=\{1\}$).  We have $\varepsilon(E)=-1$ if and only if
$E$ is unramified over $K$~; otherwise, $\varepsilon(E)>0$.  We recall
what the possibilities for $\varepsilon(E)$ are, and how it is related
to another invariant of $E$ in some cases.

\numeroter  Suppose that $e<+\infty$ and ${}_pK^\times$ has
order~$p$, and let $D\subset\overline{K^\times}$ be a line (a
 $1$-dimensional subspace).  There is a unique integer $j$ such that
 $D\subset\bar U_j$ but $D\not\subset\bar U_{j+1}$, with the
 convention that $\bar U_0=\overline{K^\times}$.  We define the {\it
 level\/} $\delta(D)$ of $D$ to be $pc-j$.  We have seen that
 $\delta(D)\in[0,pc]$, and if $\delta(D)\equiv0\pmod p$, then
 $\delta(D)=0$ or $\delta(D)=pc$.

\numeroter Suppose that $e=+\infty$ and let $D\subset\overline{K^+}$
 be a line.  There is a unique integer $j$ such that
 $D\subset\overline{\pgoth^j}$ but
 $D\not\subset\overline{\pgoth^{j+1}}$, and we define the {\it
 level\/} $\delta(D)$ of $D$ to be $-j$.  We have
 $\delta(D)\in[0,+\infty[$, and if $\delta(D)\equiv0\pmod p$, then
 $\delta(D)=0$, as we have seen.

\numeroter  {\it Suppose that\/ $e<+\infty$ and\/ ${}_pK^\times$ has
order\/~$p$ (resp.~$e=+\infty$).  Let\/ $E$ be a cyclic extension of\/ $K$ of
degree\/~$p$, and let\/  $D\subset\overline{K^\times}$
(resp.~$D\subset\overline{K^+}$) be the line such that\/ $E=K(\root
p\of D)$ (resp.~$E=K(\wp^{-1}(D))$).  If\/ $E$ is unramified over\/
$K$, then\/ $D=\bar U_{pc}$ (resp.~$D=\bar\ogoth$).  If\/ $E$ is ramified over\/
$K$, then\/ $\varepsilon(E)=\delta(D)$.} \cqfd 

\numeroter  It follows that if  $E$ is ramified over $K$, then
$\varepsilon(E)=b_p^{(i)}$ for some $i\in[1,e]$ or $\varepsilon(E)=pc$
(resp.~$\varepsilon(E)=b_p^{(i)}$ for some $i>0$), and all these
possibilites do occur.
 
\numeroter {\it Remark\/}.  The only case not covered by this
proposition is when $e<+\infty$ and ${}_pK^\times$ is trivial.  One
can compute $\varepsilon(E)$ in this case by replacing $K$ by
$K'=K(\!\root p\of1)$ and $E$ by $E'=EK'$.  If $E$ is ramified over
$K$, then $\varepsilon(E)=b_p^{(i)}$ for some $i\in[1,e]$, and all
these possibilites do occur.  In particular, $\varepsilon(E)\not\equiv0\pmod
p$, as in the case $e=+\infty$.  See \citer\locdisc(Prop.~63), for
example.  In all three cases, there are only finitely many $E$ with a
given $\varepsilon(E)$, and their number  can be easily computed.

\bigbreak
{\bf 5.  Orthogonality}
\bigskip

\numeroter Recall that for every galoisian extension $M$ of $K$, the
profinite group $G=\Gal(M|K)$ comes equipped with a decreasing
filtration $(G^t)_{t\in[-1,+\infty[}$ --- the ramification filtration
in the upper numbering --- by closed normal subgroups which is
separated ($\cap_tG^t=\{1\}$) and exhaustive ($G^{-1}=G$).  We put
$G^{t+}=\cup_{s>t}G^s$, and call $t$ a ramification break for $G$ if
$G^{t+}\neq G^{t}$, as in the case of degree-$p$ cyclic extensions
above.  In general, $G^0$ is the inertia subgroup of $G$ and $G^{0+}$
is the (wild) ramification subgroup.

\numeroter  We take $M$ to be the maximal abelian extension of
$K$ of exponent~$p$ and determine the ramification breaks of $G$.
When $e<+\infty$ and ${}_pK^\times$ has order~$p$ (resp.~$e=+\infty$),
so that $M=K(\root p\of{K^\times})$ (resp.~$M=K(\wp^{-1}(K))$), we
have the pairing
$$
G\times\overline{K^\times}\to{}_pK^\times\ (\hbox{resp.\ }
G\times\overline{K^+}\to\F_p),
$$
 as recalled in the Introduction, and we show that the filtration on
$\overline{K^\times}$ (resp.~$\overline{K^+}$) is orthogonal to the
filtration on $G$ in a certain precise sense.

\numeroter  The maximal tamely ramified extension $M_1$ of $K$ in $M$ is the
unramified degree-$p$ extension $M_1=K(\!\root p\of{U_{pc}})$
(resp.~$M_1=K(\wp^{-1}(\ogoth))$), as we have recalled.

\numeroter  {\it Suppose that\/ $e<+\infty$ and\/ ${}_pK^\times$ has
order\/~$p$.  We have\/ $G^t=G^1$ for\/ $t\in\;]-1,1]$, and, for\/
$t\in\;]0,pc+1]$,
$$
G^{t\perp}=\bar U_{pc-\lceil t\rceil+1}
$$
under\/ $G\times\overline{K^\times}\to{}_pK^\times$.  The positive
ramification breaks in the filtration on\/ $G$ occur precisely at the
$b_p^{(i)}$ ($i\in[1,e]$) and at\/ $pc$.} \cqfd

\numeroter  {\it Suppose that\/ $e=+\infty$.  We have\/ $G^t=G^1$
for\/ $t\in\;]-1,1]$, and, for\/ $t>0$,
$$
G^{t\perp}=\overline{\pgoth^{-\lceil t\rceil+1}}
$$
under\/ $G\times\overline{K^+}\to\Fp$.  The positive ramification breaks in the
filtration on\/ $G$ occur precisely at the  $b_p^{(i)}$
($i>0$).} \cqfd

\numeroter {\it Remark}.  For the proofs, see \citer\locdisc()
and \citer\further() respectively.  Now suppose that $e<+\infty$ and
${}_pK^\times$ is trivial, and put $K'=K(\!\root p\of1)$,
$\Gamma=\Gal(K'|K)$ and $M'=MK'$.  Let $\omega:\Gamma\to\F_p^\times$
be the cyclotomic character giving the action of $\Gamma$ on
${}_pK^{\prime\times}$. It can be checked that the subspace $D\subset
K^{\prime\times}\!/K^{\prime\times p}$ such that $M'=K'(\!\root p\of
D)$ is precisely the $\omega$-eigenspace for the action of $\Gamma$.
Hence or otherwise, one shows that the positive ramification breaks in
the filtration on\/ $G$ occur precisely at the $b_p^{(i)}$
($i\in[1,e]$).

\numeroter Let $L$ be an abelian extension of $K$ of exponent $p$.
It follows from the foregoing and Herbrand's theorem --- the
ramification filtration in the upper numbering is compatible with the
passage to the quotient --- that the ramification breaks of
$\Gal(L|K)$ are integers.  This is a special case of the Hasse-Arf
theorem, valid for all abelian extensions of local fields.  The
advantage of the direct proof in this special case is that one can
specify which integers occur.

\bigbreak
{\bf 6. Norms}
\bigskip

\numeroter  Suppose that $e<+\infty$ and ${}_pK^\times$ has order~$p$
(resp.~$e=+\infty$), and let $i\in[0,pc+1]$ (rep.~$i\in\N$) be an
integer.  Let $L_i=K(\!\root p\of{U_{pc-i+1}})$
(resp.~$L_i=K(\wp^{-1}(\pgoth^{-i+1}))$.  Recall that $L_0=K$, and
that $L_1$ is the unramified degree-$p$ extension of $K$.  The
inductive limit of the $L_i$ (which is nothing but $L_{pc+1}$ if
$e<+\infty$ and ${}_pK^\times$ has order~$p$) is equal to the maximal
abelian extension $M$ of $K$ of exponent~$p$.

\numeroter {\it We have\/ $N_{L_i|K}(L_i^\times)=U_iK^{\times p}$ for
every\/ $i\in[0,pc+1]$ (rep.~$i\in\N$), with the convention that\/
$U_0=K^\times$.} \cqfd

\numeroter All this can presumably
be proved by studying, as in \citer\corpslocaux(Chapter~V), what the
norm maps $N_{L_i|K}:L_i^\times\to K^\times$ do to the filtrations.
It follows that
$K^\times\!/N_{L_i|K}(L_i^\times)=\overline{K^\times}/\bar U_i$ for
every $i\in[0,pc+1]$ (rep.~$i\in\N$).  When $i=1$, the surjective
valuation on $K$ induces an isomorphism $\bar
v:K^\times\!/N_{L_1|K}(L_1^\times)\to\Z/p\Z$.

\bigbreak
{\bf 7. Reciprocity}
\bigskip

\numeroter  Keep the previous notation and continue to suppose that
$e<+\infty$ and ${}_pK^\times$ has order~$p$ (resp.~$e=+\infty$) As
$L_1$ is the unramified degree-$p$ extension of $K$, the group
$\Gal(L_1|K)$ has a canonical generator $\sigma$ --- the one which
reduces to the $k$-automorphism $x\mapsto x^q$ ($q=\Card k$) of the
residue field of $L_1$.

\numeroter {\it There is a unique isomorphism\/
$\rho_1:\overline{K^\times}/\bar U_1\to\Gal(L_1|K)$ such that\/
$\rho_1(\bar\pi)=\sigma$ for every uniformiser\/ $\pi$ of\/ $K$. The
kernel of the resulting map\/ $K^\times\to\Gal(L_1|K)$ is\/
$N_{L_1|K}(L_1^\times)$}. \cqfd

\numeroter For $i\in[1,pc+1]$ (resp.~$i>0$), and for every
intermediate extensoin $K\subset E\subset L_i$ we have the galoisian
projection $\Gal(L_i|K)\to\Gal(E|K)$.  In particular, we hav the
projection $\Gal(L_i|K)\to\Gal(L_1|K)$.  Recall also that
$K^\times\!/N_{L_i|K}(L_i^\times)=\overline{K^\times}/\bar U_i$.

\numeroter  {\it For\/ $i\in[1,pc+1]$ (resp.~$i>0$), there is a unique
isomorphism\/ $\rho_i$ making the square
$$
\diagram{
\overline{K^\times}/\bar U_i&\droite{\rho_i}&\Gal(L_i|K)\cr
\vfl{}{}{5mm}&&\vfl{}{}{5mm}\cr
\overline{K^\times}/\bar U_1&\droite{\rho_1}&\Gal(L_1|K)\cr
} 
$$
commute, and such that for every intermediate extensoin\/ $K\subset
E\subset L_i$, the kernel of the map\/ $K^\times\to\Gal(E|K)$ deduced
from $\rho_i$ is\/ $N_{E|K}(E^\times)$.}  \cqfd

\numeroter Let\/ $M=K(\root p\of{K^\times})$ (resp.~$M=K(\wp^{-1}(K))$)
be the maximal abelian extension of\/ $K$ of exponent\/~$p$.  It
follows from the foregoing that there is a unique isomorphism\/
$\rho:\overline{K^\times}\to\Gal(M|K)$ of profinite groups such that
the resulting map $K^\times\to\Gal(L_1|K)$ takes every uniformiser
of\/ $K$ to\/ $\sigma$ and such that for every intermediate
extensoin\/ $K\subset E\subset M$ of finite degree, the kernel of the
resulting map\/ $K^\times\to\Gal(E|K)$ is\/ $N_{E|K}(E^\times)$.
Moreover, $\rho(\bar U_i)=\Gal(M|K)^i$ for every $i\in[1,pc+1]$
(resp. $i>0$).

\bigbreak
{\bf 8. Orthogonality bis}
\bigskip

\numeroter  Suppose that $e<+\infty$ and ${}_pK^\times$ has
order~$p$.  Combining the kummerian pairing
$G\times\overline{K^\times}\to{}_pK^\times$ with the reciprocity
isomorphism $\rho:\overline{K^\times}\to G$, we get the hilbertian
pairing $\overline{K^\times}\times\overline{K^\times}\to{}_pK^\times$.

\numeroter {\it For every\/ $i\in[0,pc+1]$, the orthogonal complement
of\/ $\bar U_i$ for the hilbertian pairing is\/ $\bar U_{pc-i+1}$.} \cqfd

\numeroter  Suppose finally that $e=+\infty$.  As before, combining the  pairing
$G\times\overline{K^+}\to\F_p$ with the reciprocity isomorphism
$\rho:\overline{K^\times}\to G$, we get the hilbertian pairing
$\overline{K^\times}\times\overline{K^+}\to\F_p$.

\numeroter {\it For every\/ $i\in\N$, the orthogonal complement
of\/ $\bar U_i$ for the hilbertian pairing is\/
$\overline{\pgoth^{-i+1}}$ and vice versa.} \cqfd

\numeroter  {\it Remark}. In both cases, $\bar U_0=\overline{K^\times}$
by convention.  This Note was written in response to a question by
Hatice Boylan on MathOverflow \citer\boylan().

\bigbreak
\unvbox\bibbox

\bye